\documentclass[12pt]{article}
\usepackage[dvips]{epsfig}
\usepackage{graphicx}
\usepackage{amsthm,amsmath,amssymb}
\usepackage{mathtools}
\usepackage{latexsym}
\usepackage{graphics}
\usepackage{url}
\usepackage{placeins}
\usepackage{color}
\usepackage{bm}
\usepackage{comment}
\usepackage{statex}
\usepackage{pgfplots}

\usepackage[
            bookmarksnumbered=true,
            bookmarksopen=true,
            colorlinks=true,
            citecolor=green,
            linkcolor=green,
            anchorcolor=red,
            urlcolor=blue
            ]{hyperref}

\definecolor{blue}{rgb}{0,0,0.9}
\definecolor{red}{rgb}{0.9,0,0}
\definecolor{green}{rgb}{0,0.50,0.10}
\definecolor{violet}{rgb}{0.5804,0.0000,0.8275}

\def\@themcountersep{}

\makeatletter
\newcommand{\labeltext}[2]{%
  \@bsphack
  \csname phantomsection\endcsname 
  \def\@currentlabel{#1}{\label{#2}}%
  \@esphack
}
\makeatother

\newtheorem{THEO}{Theorem}[section]
\newtheorem{ALGo}[THEO]{Algorithm}
\newtheorem{ASSUMPT}[THEO]{Assumption}
\newtheorem{CONJ}[THEO]{Conjecture}
\newtheorem{COND}[THEO]{Condition}
\newtheorem{CORO}[THEO]{Corollary}
\newtheorem{DEFI}[THEO]{Definition}
\newtheorem{EXAMP}[THEO]{Example}
\newtheorem{INSTANCE}[THEO]{Instance}
\newtheorem{FACT}[THEO]{Fact}
\newtheorem{HYPO}[THEO]{Hypothesis}
\newtheorem{LEMM}[THEO]{Lemma}
\newtheorem{PROB}[THEO]{Problem}
\newtheorem{PROP}[THEO]{Proposition}
\newtheorem{REMA}[THEO]{Remark}
\newcommand{\theo}{\begin{THEO}}
\newcommand{\algo}{\begin{ALGo} \rm}
\newcommand{\assumpt}{\begin{ASSUMPT} \rm}
\newcommand{\cond}{\begin{COND}}
\newcommand{\conj}{\begin{CONJ}}
\newcommand{\coro}{\begin{CORO}}
\newcommand{\defi}{\begin{DEFI} \rm}
\newcommand{\examp}{\begin{EXAMP} \rm}
\newcommand{\instan}{\begin{INSTANCE} \rm}
\newcommand{\fact}{\begin{FACT}}
\newcommand{\hypo}{\begin{HYPO} \rm}
\newcommand{\lemm}{\begin{LEMM}}
\newcommand{\prob}{\begin{PROB} \rm}
\newcommand{\prop}{\begin{PROP}}
\newcommand{\rema}{\begin{REMA} \rm}

\newcommand{\etheo}{\end{THEO}}
\newcommand{\ealgo}{\end{ALGo}}
\newcommand{\eassumpt}{\end{ASSUMPT}}
\newcommand{\econd}{\end{COND}}
\newcommand{\econj}{\end{CONJ}}
\newcommand{\ecoro}{\end{CORO}}
\newcommand{\edefi}{\end{DEFI}}
\newcommand{\eexamp}{\end{EXAMP}}
\newcommand{\einstan}{\end{INSTANCE}}
\newcommand{\efact}{\end{FACT}}
\newcommand{\ehypo}{\end{HYPO}}
\newcommand{\elemm}{\end{LEMM}}
\newcommand{\eprob}{\end{PROB}}
\newcommand{\eprop}{\end{PROP}}
\newcommand{\erema}{\end{REMA}}

\def\0{\mbox{\bf 0}}
\def\1{\mbox{\bf 1}}
\def\2{\mbox{\bf 2}}
\def\3{\mbox{\bf 3}}
\def\4{\mbox{\bf 4}}
\def\5{\mbox{\bf 5}}
\def\6{\mbox{\bf 6}}
\def\7{\mbox{\bf 7}}
\def\8{\mbox{\bf 8}}
\def\9{\mbox{\bf 9}}
\def\a{\mbox{\boldmath $a$}}

\def\u{\mbox{\boldmath $u$}}
\def\v{\mbox{\boldmath $v$}}
\def\w{\mbox{\boldmath $w$}}
\def\x{\mbox{\boldmath $x$}}
\def\y{\mbox{\boldmath $y$}}

\def\A{\mbox{\boldmath $A$}}
\def\B{\mbox{\boldmath $B$}}
\def\C{\mbox{\boldmath $C$}}

\def\O{\mbox{\boldmath $O$}}
\def\P{\mbox{\boldmath $P$}}

\def\U{\mbox{\boldmath $U$}}
\def\V{\mbox{\boldmath $V$}}
\def\W{\mbox{\boldmath $W$}}
\def\X{\mbox{\boldmath $X$}}

\def\EC{\mbox{$\cal E$}}

\def\VC{\mbox{$\cal V$}}

\def\inprod#1#2{\langle#1, \, #2\rangle}
\def\Binprod#1#2{{\Bigl\langle}#1, \, #2{\Bigr\rangle}}

\def\bgamma{\mbox{\boldmath $\gamma$}}

\def\bdelta{\mbox{\boldmath $\delta$}}

\def\s0{\mbox{\scriptsize \boldmath $0$}}

\def\Real{\mathbb{R}}

\def\SymMat{\mathbb{S}}

\begin{document}

\title{Separable QCQPs and Their Exact SDP Relaxations}

\author{
\normalsize
Masakazu Kojima\thanks{Department of Data Science for Business Innovation, Chuo University, Tokyo, Japan ({\tt kojima@is.titech.ac.jp}).} \and \normalsize
Sunyoung Kim\thanks{Department of Mathematics, Ewha W. University, 
Seoul, 
Korea 
			({\tt skim@ewha.ac.kr}). 
			 } \and \normalsize
Naohiko Arima\thanks{
	({\tt nao$\_$arima@me.com}).} 
}

\date{\today}

\maketitle 


\begin{abstract}
\noindent
This paper studies exact semidefinite programming relaxations (SDPRs) for separable quadratically constrained quadratic programs (QCQPs).
We consider the construction of a larger separable QCQP from multiple QCQPs with exact SDPRs.
We show that exactness is preserved when such QCQPs are combined through a separable horizontal connection,
where the coupling is induced through the right-hand-side parameters of the constraints.
The proposed framework provides a simple sufficient condition for exactness of the resulting SDPR.
We then identify notable classes of QCQPs for which this condition holds,
including convex QCQPs, QCQPs defined by sign-pattern and graph-structural conditions, and separable homogeneous QCQPs with a limited number of constraints.
Two examples illustrate the constructive nature of the proposed framework, showing how heterogeneous QCQPs can be combined to yield new instances with exact SDP relaxations.
\end{abstract}


\vspace{0.5cm}




\section{Introduction} 

\label{section:introduction}

Quadratically constrained quadratic programs (QCQPs) 
form a fundamental class of generally nonconvex and computationally challenging 
optimization problems (see, e.g., \cite{FLOUDAS1995}).
They arise in a wide range of applications, including signal processing \cite{LUO2010,HUANG2010}, control \cite{PAPPAS2021}, finance \cite{JORION2003}, and combinatorial optimization \cite{LOIOLA2007,RENDLE10}. 
Among convex relaxation techniques for QCQPs, semidefinite programming relaxations (SDPRs) have been studied extensively as both a theoretical framework and a practical tool for computing lower bounds and approximate solutions. 
 In general, the optimal value of an SDPR provides a lower bound for that of the original QCQP.
Let $\eta$ and $\zeta$ denote the optimal values of the SDPR and the original QCQP, respectively. 
The SDPR is said to be exact if $\eta = \zeta$. 
When an optimal solution of the QCQP exists, the SDPR is exact if and only if it admits a rank-one optimal 
solution corresponding to an optimal solution of the original QCQP.


Over the past two decades, 
many structured classes of QCQPs have been identified for which the standard SDPR is exact. Representative examples include 
\vspace{-2mm}
\begin{description}
\item[(a) ] QCQPs characterized by convexity \cite{BENTAL2001}, \vspace{-2mm} 
\item[(b) ] QCQPs characterized by sign-pattern conditions and associated graph-structural conditions 
\cite{AZUMA2023,SOJOUDI2014},\vspace{-2mm} 
\item[(c) ] 
Separable homogeneous QCQPs characterized by a limited number of constraints 
\cite{HUANG2010,LUO2010}, 
\vspace{-2mm} 
\item[(d) ] QCQPs characterized by the rank-one-generated (ROG) property  
\cite{ARGUE2023,KIM2020},\vspace{-2mm} 
\item[(e) ] QCQPs characterized by various non-intersecting quadratic constraint (NIQC) conditions 
\cite{ARIMA2024,JOYCE2024,YANG2018}. \vspace{-2mm} 
\end{description}
These classes have largely been studied independently. Recently, \cite{ARIMA2024}
 investigated the relationship between classes (d) and (e), showing that (e) can be regarded as a special case of (d).

In parallel with this line of research, Kojima, Kim, and Arima \cite{KOJIMA2025} proposed a unified framework for extending QCQPs 
that admit exact convex relaxations 
by adding quadratic inequality constraints under suitable NIQC assumptions, without changing the variable set. 
Starting from QCQPs in classes (a)--(e), their framework provides sufficient conditions under which the extended problem continues to
admit an exact SDP, CPP, or DNN relaxation. 
In this sense, the framework enlarges the class of QCQPs with exact convex relaxations ``vertically'' by adding constraints 
while preserving exactness.


In contrast to the above line of work, which extends QCQPs with exact SDPRs by adding constraints, 
we pursue a different but complementary direction: 
combining QCQPs with exact SDPRs through a separable coupling structure. 
Specifically, we consider $\hat{p}$ QCQPs, 
each parameterized by the right-hand-side vector of its constraints, and construct a larger QCQP by coupling these problems through shared constraint parameters.
We interpret this construction as a horizontal connection of QCQPs with exact SDPRs.


More precisely, for each $p=1,\ldots,\hat{p}$, consider the QCQP 
\begin{eqnarray}
\zeta^p(\bdelta^p) & = & \inf\left\{ f^p_0(\u^p) : f^p_k(\u^p) \trianglelefteq_k \delta^p_k \ (1 \leq k \leq m) \right\}.
\label{eq:QCQPeach}
\end{eqnarray}
Here, $f^p_k$ denotes a real-valued quadratic function on the $n^p$-dimensional  Euclidean space $\Real^{n^p}$, $\bdelta^p=(\delta^p_1,\ldots,\delta^p_m)\in\Real^m$, 
and $\trianglelefteq_k$ denotes either `$\le$', `$=$', or `$\ge$'. Without loss of generality, we assume that $f^p_k(\0)=0$ 
$(0 \leq k \leq m,1 \leq p \leq \hat{p})$. 
We also note that each QCQP may involve a redundant constraint such that $f^p_{k}$ and $\delta_k$ 
are identically $0$ for some $k \in \{1,\ldots,m\}$. Therefore,
the number of nonredundant constraints may vary across the QCQPs. The resulting horizontally connected QCQP is
\begin{eqnarray}
\zeta(\bgamma) & = & \inf\left\{ \sum_{p=1}^{\hat{p}} f^p_0(\u^p) : \sum_{p=1}^{\hat{p}} f^p_k(\u^p) \trianglelefteq_k \gamma_k \ (1 \leq k \leq m) \right\},
\label{eq:QCQP00}
\end{eqnarray}
where $\bgamma=(\gamma_1,\ldots,\gamma_m)\in\Real^m$. We call each QCQP~\eqref{eq:QCQPeach} a sub-QCQP of QCQP~\eqref{eq:QCQP00}, 
and refer to this construction as {\it a horizontal connection of sub-QCQPs}  (hereafter, the horizontal connection).  This contrasts with vertical extensions, 
which enlarge QCQPs by adding constraints without changing the variable set.


A central feature of the horizontal connection is the separable structure of both the objective and 
the constraints. QCQPs with this structure arise, for example, in unicast downlink transmit beamforming 
\cite{BENGTSSON2002,HUANG2010,LUO2010}. In that setting, existing exactness results concern 
the homogeneous case, where all functions $f_k^p \ (0 \le k \le m, 1 \le p \le \hat{p})$ are homogeneous quadratic functions, 
and show   
that the associated SDPR is exact 
under conditions including that the number m of constraints does not exceed $\hat{p} + 1$, i.e.,
$m \le \hat{p}+1$.
This is precisely the class of separable homogeneous QCQPs with a limited number of constraints described in (c) above. 
The following theorem extends this exactness result by allowing 
more general classes of sub-QCQPs within the horizontal connection framework. The theorem shows that, under mild assumptions, exactness of the SDPR is preserved under the horizontal 
connection of sub-QCQPs whose SDPRs are exact.


\theo \label{theorem:main0}
Assume that \vspace{-2mm} 
\begin{description}
\item{(I) } for each $p = 1,\ldots,\hat{p}$, the SDPR~\eqref{eq:SDPRp} of QCQP~\eqref{eq:QCQPeach} is exact 
whenever $-\infty < \eta^p(\bdelta^p) < \infty$, and \vspace{-2mm} 
\item{(II)} the SDPR~\eqref{eq:SDPR0}  of QCQP~\eqref{eq:QCQP00} admits an optimal solution.  \vspace{-2mm}
\end{description}
Then the SDPR~\eqref{eq:SDPR0} is exact. 
\etheo
\noindent
Assumption (I) is weaker  than requiring exactness for every $\bdelta^p \in \Real^m$, since it only requires exactness 
when the optimal value of SDPR~\eqref{eq:SDPRp} is finite. This distinction is essential, as the SDPR may attain 
the value $-\infty$ even when the original QCQP has a finite optimal value. Such a phenomenon can occur even in simple 
instances; see \cite[Example 4.2]{KIM2023}.


Theorem~\ref{theorem:main0} applies not only to class (c) but also to QCQPs in classes (a) and (b). Indeed, 
for these classes, the exactness of the SDPR depends only on the quadratic and linear terms and 
is independent of the constant terms in the quadratic inequality constraints. 
This property is consistent with the horizontal connection framework, where sub-QCQPs 
are coupled through their right-hand-side parameters. 
Thus, Theorem 1.1 provides a unified way to construct separable QCQPs with exact SDPR 
by combining heterogeneous sub-QCQPs from classes (a), (b), and (c). 
In contrast, classes (d) and (e) do not fit this framework directly,  
since their defining conditions depend essentially on the constant terms of the quadratic inequalities.


QCQP~\eqref{eq:QCQP00} also admits an interpretation as a bilevel optimization problem:
\begin{eqnarray}
\left.
\begin{array}{ll}
\mbox{Minimize }  &  \zeta^1(\v^1)+ \cdots +  \zeta^{\hat{p}}(\v^{\hat{p}}) \\ 
\mbox{subject to } &  \v^1+ \cdots + \v^{\hat{p}}= \bgamma, \\
                              &  \zeta^p(\v^p) 
                              = \inf\left\{ f^p_{0}(\u^p) : f^p_{k}(\u^p) \trianglelefteq_k v^p_k \ (1 \leq k \leq m) \right\} 
                              \ (1 \leq p \leq \hat{p}). 
 \end{array}  
 \right\} \label{eq:bilevel}                        
\end{eqnarray}
Here, the upper-level problem allocates the parameters $\v^1,\ldots, \v^{\hat{p}}$ subject to 
$\v^1+\cdots+ \v^{\hat{p}}  = \bgamma$, while each lower-level problem computes the optimal 
value of a QCQP whose constraints are parameterized by the corresponding 
right-hand-side vector $\v^p$. 
Thus, QCQP~\eqref{eq:QCQP00} may be viewed as a special class of 
bilevel optimization problems. 
In contrast to standard resource allocation problems~\cite{BERTSEKAS1999,LUO2010}, 
the variables $\v^p$ are not restricted to be nonnegative. 
This difference simplifies the problem structure, as it can be reduced to 
the single-level separable QCQP~\eqref{eq:QCQP00}.
Theorem~\ref{theorem:main0} ensures that, under the stated assumptions, the corresponding 
single-level formulation~QCQP~\eqref{eq:QCQP00} admits an exact SDPR.

\rema 
In Theorem~\ref{theorem:main0}, we have assumed neither the existence of a rank-one optimal solution for 
SDPR~\eqref{eq:SDPRp}  in Assumption I nor that for 
SDPR~\eqref{eq:SDPR0} in Assumption II. 
Under our definition, exactness refers only to the equality of the optimal values of a QCQP and its SDPR, and does not 
guarantee the existence of a rank-one optimal solution of the SDPR 
corresponding to an optimal solution of the QCQP in general. 
See \cite[Example 4.3]{KIM2023}. 
\erema


The main contribution of this paper is to show that the exactness of the SDPR is preserved under a horizontal connection of QCQPs whose individual SDPRs are exact, 
thereby providing a unified framework for constructing separable QCQPs with exact SDP relaxations from heterogeneous building blocks.
In addition, we establish a rank-based exactness result for a class of separable homogeneous QCQPs, extending existing results based on constraint bounds \cite{LUO2010}.
Although the present paper does not focus on specific applications, the proposed horizontal connection framework provides a systematic mechanism for constructing 
larger-scale QCQPs with guaranteed exact SDP relaxations from smaller building blocks. Such constructions are potentially useful in applications where complex 
systems are naturally decomposed into loosely coupled subsystems.
To complement the theoretical development, we present two examples that demonstrate how the proposed framework can be used constructively. 
The first example examines a separable homogeneous QCQP, while the second illustrates how heterogeneous sub-QCQPs can be combined 
to form a larger QCQP with an exact SDP relaxation. These examples highlight the role of the framework as a tool for generating new problem instances with guaranteed exactness.

%

The remainder of the paper is organized as follows. 
Section~\ref{section:SDPR} introduces the SDPR formulations and notation,
Section~\ref{section:proof} proves Theorem~\ref{theorem:main0},
Section~\ref{section:subQCQPs} discusses representative classes of sub-QCQPs satisfying Assumption (I),
Section~\ref{section:example} presents two examples illustrating the main results, and
Section~\ref{section:ConcludingRemarks} concludes the paper.




\section{SDPRs of QCQPs~\eqref{eq:QCQPeach} and~\eqref{eq:QCQP00}}

\label{section:SDPR}

SDPR~\eqref{eq:SDPRp} is the standard Shor semidefinite programming relaxation of QCQP~\eqref{eq:QCQPeach}. 
If this relaxation admits a rank-one optimal solution, then it 
yields an optimal solution of the original QCQP. 
For QCQP~\eqref{eq:QCQP00}, 
the separable structure induces a chordal sparsity structure 
in the Shor relaxation. 
Exploiting this structure, the lifted 
matrix variable can be decomposed into block matrices, leading to the sparse 
SDPR~\eqref{eq:SDPR0}. In this formulation, a block-wise rank-one optimal solution corresponds 
to an optimal solution of QCQP~\eqref{eq:QCQP00}. 
This construction can be viewed as a special case of the matrix completion framework of~\cite{FUKUDA2003}; 
see also Remark~\ref{remark:blockdiagonal}.


To describe the SDPRs of QCQPs~\eqref{eq:QCQPeach} and~\eqref{eq:QCQP00}, we  first introduce  notation. 
Let $\SymMat^n$ denote the space of $n \times n$ symmetric matrices equipped with the inner product
\[
\inprod{\A}{\B} = \sum_{i=1}^n \sum_{j=1}^n A_{ij} B_{ij} \ \text{for } \A, \B \in \SymMat^n,
\]
and let $\SymMat^n_+$ denote the cone of $n \times n$ symmetric positive semidefinite matrices.
We represent each quadratic function $f^p_k : \Real^{n^p} \rightarrow \Real$ 
$(1 \leq p \le \hat{p},0\le k \le m)$ as 
\[
f_k^p(\u^p) =
\begin{pmatrix} \u^p \\ 1 \end{pmatrix}^T
\B_k^p
\begin{pmatrix} \u^p \\ 1 \end{pmatrix}
\quad \text{for every } \u^p \in \Real^{n^p},
\]
where $\B_k^p \in \SymMat^{n^{p+1}}$ and $[B_k^p]_{n^p+1,n^p+1} = 0$, reflecting the assumption $f_k^p(0) = 0$.
Equivalently, this can be expressed in inner product form as
\[
f_k^p(\u^p) =
\Binprod{\B_k^p}
{
\begin{pmatrix}
\u^p (\u^p)^T & \u^p \\
(\u^p)^T & 1
\end{pmatrix}
}
\quad \text{for every } \u^p \in \Real^{n^p}. 
\]
With this representation, QCQPs~\eqref{eq:QCQPeach} and~\eqref{eq:QCQP00}
admit the following equivalent vector and lifted-matrix formulations: 
\begin{eqnarray}
\zeta^p(\bdelta^p) 
& = & \inf\left\{ f^p_0(\u^p) : \u^p \in \Real^{n^p}, \ f^p_k(\u^p) \trianglelefteq_k \delta^p_k  \ (1 \leq k \leq m) \right\} \nonumber \\
& = & 
\inf\left\{ \Binprod{\B^p_{0}}{\begin{pmatrix} \u^p(\u^p)^T & \u^p \\ (\u^p)^T & 1 \end{pmatrix}}: 
\begin{array}{l}
\u^p \in \Real^{n^p}, \\[3pt]
 \Binprod{\B^p_{k}}{\begin{pmatrix} \u^p(\u^p)^T & \u^p \\ (\u^p)^T & 1 \end{pmatrix}} 
 \trianglelefteq_k \delta^p_k  \\
  (1 \leq k \leq m)
 \end{array}
\right\}, \label{eq:QCQPp}
\end{eqnarray}
$(1 \leq p\leq \hat{p})$ and 
\begin{eqnarray}
\zeta(\bgamma) 
& = & \inf\left\{ \sum_{p=1}^{\hat{p}} f^p_0(\u^p) : 
\begin{array}{l}
\u^p \in \Real^{n^p} \ (1 \leq p\leq \hat{p}), \\[3pt]
\displaystyle 
 \sum_{p=1}^{\hat{p}} f^p_k(\u^p) \trianglelefteq_k \gamma_k \ (1 \leq k \leq m) 
\end{array}
\right\} \nonumber \\ 
& = & \inf\left\{ \sum_{p=1}^{\hat{p}}  \Binprod{\B^p_{0}}{\begin{pmatrix} \u^p(\u^p)^T & \u^p \\ (\u^p)^T & 1 \end{pmatrix}}: 
\begin{array}{l}
\u^p \in \Real^{n^p} \ (1 \leq p\leq \hat{p}), \\[3pt]
\displaystyle \sum_{p=1}^{\hat{p}} \Binprod{\B^p_{k}}{\begin{pmatrix} \u^p(\u^p)^T & \u^p \\ (\u^p)^T & 1 \end{pmatrix}} \trianglelefteq_k \gamma_k \\ 
 (1 \leq k \leq m) 
\end{array}
\right\}. \label{eq:QCQPseparable2}
\end{eqnarray}


By replacing the rank-one matrix 
$
\begin{pmatrix}
\u^p (\u^p)^T & \u^p \\
(\u^p)^T & 1
\end{pmatrix}
$ 
with a matrix variable $\X^p \in \SymMat^{n^p+1}_+$ satisfying $X^p_{n^p+1,n^p+1} = 1$, 
we obtain the following (Shor-type) SDPRs of QCQPs~\eqref{eq:QCQPeach} and~\eqref{eq:QCQP00}: 
\begin{eqnarray}
\eta^p(\bdelta^p) 
& = & \inf\left\{   \inprod{\B^p_{0}}{\X^p}: 
\begin{array}{l}
\X^p \in \SymMat^{n^p+1}_+, \  
X^p_{n^p+1,n^p+1} = 1, \\[3pt]
 \inprod{\B^p_{k}}{\X^p} \trianglelefteq_k \delta^p_k \  (1 \leq k \leq m) 
\end{array}
\right\} 
 \label{eq:SDPRp}
\end{eqnarray}
$(1 \leq p\leq \hat{p})$, and 
\begin{eqnarray}
\eta(\bgamma) 
& = & \inf\left\{ \sum_{p=1}^{\hat{p}}  \inprod{\B^p_{0}}{\X^p}: 
\begin{array}{l}
\X^p \in \SymMat^{n^p+1}_+ \ (1 \leq p\leq \hat{p}), \\[3pt] 
X^p_{n^p+1,n^p+1} = 1  \ (1 \leq p\leq \hat{p}), \\[3pt]
\displaystyle \sum_{p=1}^{\hat{p}} \inprod{\B^p_{k}}{\X^p} \trianglelefteq_k \gamma_k \  (1 \leq k \leq m) \\
\end{array}
\right\}.  
 \label{eq:SDPR0}
\end{eqnarray}
In this relaxation, the rank-one constraint on each matrix variable is dropped.


\rema \label{remark:blockdiagonal}
This remark clarifies that the separable structure assumed in QCQP~\eqref{eq:QCQP00} 
is not intrinsic, but depends on the representation of the problem. Consider a general QCQP
\begin{eqnarray}
\zeta(\bgamma)=\inf\{f_0(\x):\x\in\Real^n,\ f_k(\x)\le \gamma_k\ (1\le k\le m)\}, 
\label{eq:QCQPgeneral}
\end{eqnarray}
where
\[
f_k(\x)=\x^T \A_k\x+\a_k^T \x
\quad \text{for every } x\in\Real^n,
\quad \A_k\in S^n,\ \a_k\in\Real^n
\quad (0\le k\le m).
\]
It is straightforward to see that QCQP~\eqref{eq:QCQPgeneral} is invariant under orthogonal linear transformations. 
Namely, for any $n \times n$ orthogonal matrix $\P$, QCQP~\eqref{eq:QCQPgeneral} is equivalent to
\begin{eqnarray}
\zeta'(\bgamma)=\inf\{f_0(\P\y):\y\in\Real^n,\ f_k(\P\y)\le \gamma_k\ (1\le k\le m)\}.
\label{eq:QCQPgeneral2}
\end{eqnarray}
Thus, a given QCQP may admit different equivalent representations. In particular, 
whether and how one can choose an $n \times n$ orthogonal matrix $\P$ so that QCQP~\eqref{eq:QCQPgeneral2} 
exhibits a separable structure is an interesting problem. This problem reduces to finding an $n \times n$ 
orthogonal matrix $\P$ such that
$ 
\P^T \A_0 \P,\dots,\P^T\A_m \P
$ 
share a common block-diagonal structure. 
Such a block-diagonal structure of the matrices directly induces a separable structure of QCQP~\eqref{eq:QCQPgeneral2}, 
in which each block corresponds to a sub-QAP.
See \cite{MUROTA2010} and the references therein for more details. 
\erema


\section{Proof of Theorem~\ref{theorem:main0}}

\label{section:proof}

Let $(\widetilde \X^1,\ldots,\widetilde{\X}^{\hat{p}})$ be an optimal solution of SDPR~\eqref{eq:SDPR0}, 
the existence of which  is guaranteed by Assumption (II).
For every $p \in \{1,\ldots,\hat{p}\}$, let
\begin{eqnarray*}
& & \bdelta^p = (\delta^p_1,\ldots,\delta^p_m), \
 \delta^p_k = \inprod{\B^{p}_{k}}{\widetilde{\X}^{p}} 
\ (1 \leq k \leq m). 
\end{eqnarray*}
By construction, each $\widetilde{\X}^p$ satisfies the constraints of SDPR~\eqref{eq:SDPRp}  with the right-hand side $\bdelta^p$; 
hence each $\widetilde{\X}^p$  is a feasible solution of SDPR~\eqref{eq:SDPRp}. 


We claim that $\widetilde{\X}^p$ is in fact optimal for SDPR~\eqref{eq:SDPRp} with the right-hand side $\bdelta^p$. Suppose not. 
Then there exists a feasible solution $\overline{\X}^p \in \SymMat^{n^p+1}_+$ of SDPR~\eqref{eq:SDPRp} such that
$ 
\inprod{\B_0^p}{\overline{\X}^p} < \inprod{\B_0^p}{\widetilde{\X}^p} 
$ 
and 
$ 
\inprod{\B_k^p}{\overline{\X}^p} \trianglelefteq_k \delta_k^p
\quad (1 \leq k \leq m). 
$ 
Hence 
\begin{eqnarray*}
& & \inprod{\B_0^p}{\overline{\X}^p}+ \sum_{p' \not= p} \inprod{\B_0^{p'}}{\widetilde{\X}^{{p'}} }
< \sum_{p' = 1}^{\hat{p}} \inprod{\B_0^{{p'}}}{\widetilde{\X}^{{p'}}} = \eta(\bgamma), \\
& & \inprod{\B_k^p}{\overline{\X}^p} + 
\sum_{p' \not= p} \inprod{\B_k^{{p'}}}{\widetilde{\X}^{p'}}  \trianglelefteq_k 
\delta^p_k +  \sum_{p' \not= p} \delta^{{p'}}_k \trianglelefteq_k  \gamma_k 
\quad (1 \leq k \leq m). 
\end{eqnarray*}
Thus, replacing $\widetilde{\X}^p$ in the optimal solution $(\widetilde{\X}^1,\ldots,\widetilde{\X}^{\hat{p}})$ of SDPR~\eqref{eq:SDPR0}
by $\overline{\X}^p$,  %
we obtain a feasible solution of SDPR~\eqref{eq:SDPR0} with a strictly smaller objective value than the optimal value 
$\eta(\bgamma)$,  a contradiction. 
Therefore we have shown that 
$ 
-\infty < \eta^p(\bdelta^p)=\inprod{\B_0^p}{\widetilde{\X}^p} < \infty
\quad (1 \leq p \leq \hat{p}).
$ 


By Assumption (I),  $\eta^p(\bdelta^p) = \zeta^p(\bdelta^p)$ $(1 \leq p \leq \hat{p})$. Thus  
\begin{eqnarray*}
\eta(\bgamma) = \sum_{p=1}^{\hat{p}} \inprod{\B^p_{0}}{\widetilde{\X}^p} = \sum_{p=1}^{\hat{p}} \eta^p(\bdelta^p) = \sum_{p=1}^{\hat{p}} 
\zeta^p(\bdelta^p). 
\end{eqnarray*}
Since any collection $\{\u^p: 1 \leq p \leq \hat{p}\}$ of feasible solutions of  QCQP~\eqref{eq:QCQPp} 
yields a feasible solution $(\u^1,\ldots,\u^{\hat{p}})$ of QCQP~\eqref{eq:QCQPseparable2},  
$\eta(\bgamma) = \sum_{p=1}^{\hat{p}} \zeta^p(\bdelta^p) \geq \zeta(\bgamma)$. 
We also know that $\eta(\bgamma) \leq \zeta(\bgamma)$ in general. Therefore $\eta(\bgamma) = \zeta(\bgamma)$, {\it i.e.}, 
SDPR~\eqref{eq:SDPR0} is exact. 
\qed




\section{Eligible classes of sub-QCQPs with exact  SDPRs}

\label{section:subQCQPs}

Section~\ref{section:introduction} describes three classes of QCQPs relevant to Theorem~\ref{theorem:main0}: (a) convex QCQPs, 
(b) QCQPs characterized by sign-pattern and graph-structural conditions, and (c) separable homogeneous QCQPs 
with a limited number of constraints. 
In this section, we show that these classes can serve as eligible sub-QCQPs in the framework of Theorem~\ref{theorem:main0}. 
Specifically, Sections~4.1, 4.2, and 4.3 discuss these three classes, respectively. Whenever the sub-QCQPs belong to any of these 
classes, Assumption~(I) of Theorem~\ref{theorem:main0} is satisfied; if Assumption~(II) also holds, then SDPR~\eqref{eq:SDPR0} 
is exact, {\it i.e.}, $\eta(\bgamma)=\zeta(\bgamma)$.


For simplicity of notation, we drop the superscript $p$ and write a sub-QCQP embedded in QCQP~\eqref{eq:QCQP00} as
\begin{eqnarray}
\zeta(\bdelta) 
& = & \inf\left\{ f_{0}(\u) : \u \in \Real^{n}, \ f_{k}(\u) \trianglelefteq_k \delta_k  \ (1 \leq k \leq m) \right\} \nonumber \\
& = & 
\inf\left\{ \inprod{\B_{0}}{\begin{pmatrix} \u\u^T & \u \\ \u^T & 1 \end{pmatrix}}: 
\begin{array}{l}
\u \in \Real^{n}, \\[3pt]
 \inprod{\B_{k}}{\begin{pmatrix} \u\u^T & \u \\ \u^T & 1 \end{pmatrix}} \trianglelefteq_k \delta_k  \ (1 \leq k \leq m)
 \end{array}
\right\}, \label{eq:QCQPp2}
\end{eqnarray}
and its SDPR  as 
\begin{eqnarray}
\eta(\bdelta) 
& = & \inf\left\{   \inprod{\B_{0}}{\X}: 
\begin{array}{l}
\X \in \SymMat^{n+1}_+, \  [\X]_{n+1,n+1} = 1, \\[3pt]
 \inprod{\B_{k}}{\X} \trianglelefteq_k \delta_k \  (1 \leq k \leq m) 
\end{array}
\right\}.   \label{eq:SDPRp2}
\end{eqnarray}
Here 
\begin{eqnarray*}
& & \B_k \in \SymMat^{n+1} \ (0 \leq k \leq m),   \\
& & f_k(\u) =  \Binprod{\B_{k}}{\begin{pmatrix} \u\u^T & \u \\ \u^T & 1 \end{pmatrix}} \ 
\mbox{for every } \u \in \Real^n \ (0 \leq k \leq m). 
\end{eqnarray*}
We now present three classes of QCQP~\eqref{eq:QCQPp2} for which  SDPR~\eqref{eq:SDPRp2} is exact.

\subsection{QCQPs characterized by convexity} 

\label{section:convexity}


We assume that `$\trianglelefteq_k$' $=$ `$\le$' $(1 \leq k \leq m)$, and 
that the objective and constraint functions are convex. 
Then each QCQP~\eqref{eq:QCQPp2} is a convex quadratic optimization problem.  
For this class of QCQPs, SDPR~\eqref{eq:SDPRp2} is exact for 
every right-hand-side vector $\bdelta \in \Real^m$, 
and hence Assumption~(I) of Theorem~\ref{theorem:main0} is automatically satisfied.
The following result is well-known (\cite{SHOR1987},\cite[ Section 4.2]{BENTAL2001}), 
 and can be proved easily. 
 \theo \label{theorem:convex} 
 Assume that $f_k : \Real^n \rightarrow \Real$ is convex $(0 \leq k \leq m)$. \vspace{-1mm} 
 \begin{description}
 \item{(i) }
 Let $\overline{\X} = \begin{pmatrix} \overline{\U} & \bar{\u} \\ \bar{\u}^T & 1 \end{pmatrix}$ 
 be an optimal solution of SDPR \eqref{eq:SDPRp2}. Then $\bar{\u}$ is an optimal solution 
 of QCQP \eqref{eq:QCQPp2}.\vspace{-2mm} 
 \item{(ii) } Let $\bar{\u} \in \Real^{n}$ be an optimal solution of QCQP \eqref{eq:QCQPp2}. 
 Then $\overline{\X} = \begin{pmatrix} \bar{\u}\bar{\u}^T & \bar{\u} \\ \bar{\u}^T & 1 \end{pmatrix} \in \SymMat^{n+1}_+$ 
 is an optimal solution of SDPR \eqref{eq:SDPRp2}. 
 \end{description}
 \etheo
 \noindent
 In particular, for this class, the exactness of the SDPR  is independent of  
 the right-hand-side vector $\bdelta$. 
 
\subsection{QCQPs characterized by sign-pattern conditions and associated graph structures}

\label{section:signPattern}

We next  consider a class of generally nonconvex QCQPs 
whose SDPR exactness is guaranteed by sign-pattern and graph-structural conditions on the coefficient matrices 
$\B_k$ $(0\leq k \leq m)$ 
and the associated sparsity 
graph. 
We assume that $\trianglelefteq_k=$ `$\le$' $(1 \leq k \leq m)$. 
As in the convex case, the exactness result for this class holds for every right-hand-side vector 
$\bdelta$. Hence these QCQPs also satisfy Assumption~(I)  of Theorem~\ref{theorem:main0}.


Let $\VC = \{1,\ldots,n\}$ and 
$ 
\EC = \{ (i,j) \in \VC \times \VC: i \not= j , \ [\B_k]_{ij} \not= 0 \ \mbox{for some } k \in \{0,1,\ldots,m\}\}.
$ 
We call $G(\VC,\EC)$ {\em the aggregated sparsity pattern graph} associated with   $\B_k$ $(0 \leq k \leq m)$.
For every $(i,j) \in \EC$, define 
\begin{eqnarray*}
\sigma_{ij} =
\begin{cases}
+1 & \text{if } [\B_k]_{ij} \ge 0 \text{ for all } k \in \{0,1,\ldots,m\},\\
-1 & \text{if } [\B_k]_{ij} \le 0 \text{ for all } k \in \{0,1,\ldots,m\},\\
0  & \text{otherwise}.
\end{cases}
\end{eqnarray*}
Let $\{C_1,\ldots,C_r\}$ denote a cycle basis for $G(\VC,\EC)$. 


\theo \label{theorem:signDefinite} (\cite[Theorem 2]{SOJOUDI2014}) Assume that \vspace{-2mm} 
\begin{description}
\item{(i) } $\sigma_{ij} \in \{-1,1\} \ \mbox{for every } (i,j) \in \EC$, \vspace{-1mm} 
\item{(ii) } $\prod_{(i,j)\in C_s} \sigma_{ij} = (-1)^{|C_s|} \mbox{for every } s=1,\ldots,r$. \vspace{-1mm} 
\end{description}
Then $\eta(\bdelta) = \zeta(\bdelta)$ for every $\bdelta \in \Real^m$. 
\etheo

As special cases, we have the following result.

\coro \label{corollary:signDefinite} (\cite[Corollary 1]{SOJOUDI2014}) 
$\eta(\bdelta) = \zeta(\bdelta)$ for every $\bdelta \in \Real^m$ 
if one of the following holds: \vspace{-2mm}
\begin{description}
\item{(i) } The graph $G(\VC,\EC)$ is arbitrary and $\sigma_{ij} = -1$ 
for every $(i,j) \in \EC$ (or equivalently all off-diagonal elements of $\B_k$ $(0 \leq k \leq m)$ are nonpositive).\vspace{-2mm}
\item{(ii) } The graph $G(\VC,\EC)$ is forest and $\sigma_{ij} \in \{ -1,1\}$ 
for every $(i,j) \in \EC$.\vspace{-2mm}
\item{(iii) } The graph $G(\VC,\EC)$ is bipartite and $\sigma_{ij} = 1$ 
for every $(i,j) \in \EC$.\vspace{-2mm}
\end{description}
\ecoro

\subsection{Separable homogeneous QCQPs characterized by a limited number of constraints}

\label{section:separable}

We now consider a homogeneous QCQP whose objective and constraint functions share a 
separable structure. Unlike the previous two classes,  
exactness of the SDPR for this class is established through a rank
argument 
that depends on the number of constraints relative to the number of separable blocks. 
Let $\bdelta \in \Real^m$,  and for each $q=1,\ldots,\hat{q}$ and $k=0,\ldots,m$, 
$g^q_k : \Real^{\bar{n}^q} \rightarrow \Real$ be a homogeneous 
quadratic function ({\it i.e}, without linear and constant terms) such that 
\begin{eqnarray*}
g^q_k(\v^q) & = & \inprod{\C^q_k}{\v^q(\v^q)^T} \ \mbox{for every } \v^q \in \Real^{\bar{n}^q}, 
\end{eqnarray*}
where $C^q_k \in \SymMat^{\bar{n}^q}$. 
Then we consider the following QCQP and its SDPR: 
\begin{eqnarray}
\zeta(\bdelta) 
& = & 
\inf\left\{ \sum_{q=1}^{\hat{q}} g^q_0(\v^q):  
\begin{array}{l} 
\v^q \in \Real^{\bar{n}^q} \ (1 \leq q \leq \hat{q}), \\
\displaystyle \sum_{q=1}^{\hat{q}} g^q_k(\v^q)  \trianglelefteq_k \delta_k \
(1 \leq k \leq m) 
\end{array}
\right\}\nonumber  \\
& = & \inf\left\{ \sum_{q=1}^{\hat{q}} \inprod{\C^q_0}{\v^q(\v^q)^T}:  
\begin{array}{l} 
\v^q \in \Real^{\bar{n}^q} \ (1 \leq q \leq \hat{q}), \\
\displaystyle \sum_{q=1}^{\hat{q}} \inprod{\C^q_k}{\v^q(\v^q)^T}  \trianglelefteq_k \delta_k \
(1 \leq k \leq m) 
\end{array}
\right\}, \label{eq:QCQP3}
\end{eqnarray}
and 
\begin{eqnarray}
\eta(\bdelta) & = & \inf\left\{ \sum_{q=1}^{\hat{q}} \inprod{\C^q_0}{\V^q}:  
\begin{array}{l} 
\V^q \in \SymMat^{\bar{n}^q}_+ \ (1 \leq q \leq {\hat{q}}),  \\ 
\displaystyle \sum_{q=1}^{\hat{q}} \inprod{\C^q_k}{\V^q} \trianglelefteq_k \delta_k \ (1 \leq k \leq m) 
\end{array}
\right\}, \label{eq:SDPR3}
\end{eqnarray}
respectively. 

\theo \label{theorem:constLimited} Assume that 
\vspace{-2mm}
\begin{description}
\item{(A) } for every optimal solution 
$(\V^1,\ldots,\V^{\hat{q}})$ of SDPR~\eqref{eq:SDPR3}, at least $m-1$ of 
\begin{eqnarray*}
\V^q\in \SymMat^n \ (1 \leq q \leq {\hat{q}}) \  \mbox{ and }  \delta_k - \sum_{q=1}^{\hat{q}} \inprod{\C^q_k}{\V^q} \in \Real 
\ (1 \leq k \leq m) 
\end{eqnarray*}
are nonzero.  \vspace{-2mm}
\item{(B) } SDPR~\eqref{eq:SDPR3} has an optimal solution, \vspace{-2mm}
\end{description}
Then SDPR~\eqref{eq:SDPR3} admits an optimal solution 
$(\widetilde{\V}^1,\ldots,\widetilde{\V}^{\hat{q}})$ such that rank$(\widetilde{\V}^q) \leq 1$ \
 $(1 \leq q \leq {\hat{q}})$. 
Consequently, $\widetilde{\V}^q=\bar{\v}^q(\bar{\v}^q)^T$ for some $\bar{\v}^q\in\Real^{\bar n^q}$, and
$(\bar{\v}^1,\ldots,\bar{\v}^{\hat q})$ is an optimal solution of QCQP~\eqref{eq:QCQP3}.
\etheo 
\proof{Introduce slack variables $s_k \in \Real$ $(1 \leq k \leq m)$ and rewrite SDPR~\eqref{eq:SDPR3} as 
\begin{eqnarray}
\eta(\bdelta) & = & \inf\left\{ \sum_{q=1}^{\hat{q}} \inprod{\C^q_0}{\V^q}:  
\begin{array}{l} 
\V^q \in \SymMat^{\bar{n}^q}_+ \ (1 \leq q \leq {\hat{q}}),  \
0 \trianglelefteq_k s_k, \\
\displaystyle \sum_{q=1}^{\hat{q}} \inprod{\C^q_k}{\V^q} - s_k = \delta_k \ (1 \leq k \leq m) 
\end{array}
\right\}. \label{eq:SDPR4}
\end{eqnarray}
By Assumption~(B), this problem admits an optimal solution.
Then we can apply 
Theorem 2.2 of \cite{PATAKI1998} to SDPR~\eqref{eq:SDPR4}, 
which provides a rank bound for feasible solutions of an SDP, for the existence of 
an optimal solution
$(\widetilde{\V}^1,\ldots,\widetilde{\V}^{\hat{q}},\bar{s}_1,\ldots,\bar{s}_m)$ 
satisfying 
\begin{eqnarray}
\sum_{q \in Q} \frac{\mbox{rank}(\widetilde{\V}^q)(\mbox{rank}(\widetilde{\V}^q)+1)}{2} 
+ \sum_{k \in K} 1 \leq m, 
\label{eq:rankCond}
\end{eqnarray}
where $Q = \{q : \widetilde{\V}^q \ne \O\}$ and 
$K = \{ k : \bar{s}_k \ne 0\}$. 
Suppose, to the contrary,  that $\mbox{rank}(\widetilde{\V}^{\bar{q}}) \geq 2$ for some $\bar{q}\in Q$. 
Then 
\[
\displaystyle \frac{\mbox{rank}(\widetilde{\V}^{\bar{q}})(\mbox{rank}(\widetilde{\V}^{\bar{q}})+1)}{2} \geq 3. 
\]
By Assumption~(A), at least $m-1$ of  the matrices 
$\widetilde{\V}^q$ $(1\le q\le \hat q)$ and  the residuals
$
\delta_k-\sum_{q=1}^{\hat q}\inprod{\C_k^q}{\widetilde{\V}^q}
 \ (1\le k\le m)
$
are nonzero. Hence $|Q \backslash\{ \tilde{q} \} |+|K|\ge m-2$. 
Therefore,
\begin{eqnarray*}
m + 1 & = & 3 + (m-2) \\
& \leq & \frac{\mbox{rank}(\widetilde{\V}^{\bar{q}})(\mbox{rank}(\widetilde{\V}^{\bar{q}})+1)}{2} + 
\sum_{q\in Q\backslash\{{\bar{q}}\}} \frac{\mbox{rank}(\widetilde{\V}^q)(\mbox{rank}(\widetilde{\V}^q)+1)}{2} 
+ \sum_{k \in K} 1. 
\end{eqnarray*}
This contradicts~\eqref{eq:rankCond}. Thus $\mbox{rank}(\widetilde{\V}^q)\le 1$ for every $q=1,\ldots,\hat q$.


Since each $\widetilde{\V}^q\in\SymMat^{\bar n^q}_+$ has rank at most one, there exists $\bar{\v}^q\in\Real^{\bar n^q}$ such that
$
\widetilde{\V}^q=\bar{\v}^q(\bar{\v}^q)^T 
$ $(1 \leq q \leq \hat{q})$. 
Substituting these factorizations into SDPR~\eqref{eq:SDPR3} shows that
$(\bar{\v}^1,\ldots,\bar{\v}^{\hat q})$ is feasible for QCQP~\eqref{eq:QCQP3} and attains the same objective value. Hence it is an optimal solution of QCQP~\eqref{eq:QCQP3}.
\qed
} 


To better understand Assumption (A), consider the case where $\ell$ of the relations 
$\trianglelefteq_k$ $(1 \leq k \leq m)$
are equalities, and the remaining ones are inequalities.
Without loss of generality, assume that $\trianglelefteq_k = $ '$=$' for $1 \leq k \leq \ell$, 
and $\trianglelefteq_k = $ '$\le$' for $\ell+1 \leq k \leq m$.
Then, for every optimal solution, the residual $\delta_k - \sum_{q=1}^{\hat{q}} \inprod{\C^q_k}{\V^q}$ 
vanishes for $1 \leq k \leq \ell$.  Hence, Assumption (A) requires that at least $m-1$ of the matrices $\V^q$ 
$(1 \leq q \leq \hat{q})$  and the remaining residuals corresponding to inequality constraints be nonzero.
This implies that $m-1 \leq \hat{q} + (m-\ell)$, {\it i.e.}, $\ell -1 \le \hat{q}$.  
In particular, when $\ell = m$, we must have $m \leq \hat{q}+1$, 
whereas when $\ell \leq 2$, Assumption (A) can be satisfied for arbitrary values of $m \geq \ell$ and $\hat{q} \geq 1$.


\rema 
Theorem~\ref{theorem:constLimited} can be compared with the result in \cite{LUO2010}, 
where the assumptions 
$m \leq {\hat{q}}+1$ and \vspace{-2mm}
\begin{description}
\item{(ii)' } $\V^q \not= \O$ $(1 \leq q \leq {\hat{q}})$ for every optimal solution $(\V^1,\ldots,\V^{\hat{q}})$ of SDPR~\eqref{eq:SDPR3} \vspace{-2mm}
\end{description}
were imposed. Under these assumptions, 
\eqref{eq:rankCond} implies 
\begin{eqnarray*}
m-1 \leq {\hat{q}} \leq \sum_{q=1}^{\hat{q}} \frac{\mbox{rank}(\widetilde{\V}^q)(\mbox{rank}(\widetilde{\V}^q)+1)}{2} 
\leq m, 
\end{eqnarray*}
and hence either $\hat q = m-1$ or $\hat q = m$; equivalently, either $m=\hat q+1$ or $m=\hat q$. 
In contrast, Assumption~(A) of Theorem~\ref{theorem:constLimited} lightens  
this restriction considerably as discussed above. 
\erema

Whether Assumption~(A) holds depends not only on the data matrices $\C_k^q$ $(1 \le q \le \hat q,\ 0 \le k \le m)$ but also 
on the right-hand-side vector $\bdelta$ and the relations `$\trianglelefteq_k$' $(1 \le k \le m)$. 
Consequently, unlike the convex and sign-pattern classes, homogeneous QCQP~\eqref{eq:QCQP3} 
does not automatically qualify as an eligible sub-QCQP in the framework of Theorem~\ref{theorem:main0}. 
When $m \leq 2$, the conclusion of Theorem~\ref{theorem:constLimited} holds without 
any additional assumption on $\bdelta$,  `$\trianglelefteq_k$' and $\C^q_k$ $(1 \leq q \leq {\hat{q}}, \ 0\leq k \leq m)$ as shown below. 
(This result is known \cite{ARIMA2024,POLYAK1998,YE2003}). Therefore, homogeneous 
QCQP~\eqref{eq:QCQP3} with $m \leq 2$ can be incorporated as an eligible sub-QCQP 
 in  separable QCQP~\eqref{eq:QCQP00} whose SDPR is guaranteed to be exact by Theorem~\ref{theorem:main0}. 

\coro \label{corollary:m=2} 
Assume that $m \leq 2$ and that SDPR~\eqref{eq:SDPR3} has an optimal solution. Then 
SDPR~\eqref{eq:SDPR3} has an optimal solution 
$(\widetilde{\V}^1,\ldots,\widetilde{\V}^{\hat{q}})$ such that rank$(\widetilde{\V}^q) \leq 1$ $(1 \leq q \leq {\hat{q}})$. \vspace{-2mm} 
\ecoro

\proof{
If Assumption~(A) of Theorem~\ref{theorem:constLimited} holds with $m=2$, then 
the conclusion follows immediately from Theorem~\ref{theorem:constLimited}. Otherwise, Assumption~(A) fails. Since $m-1=1$, 
there exists an optimal solution $(\widetilde{\V}^1,\ldots,\widetilde{\V}^{\hat q})$ of SDPR~\eqref{eq:SDPR3} 
such that all matrices $\widetilde{\V}^q$ $(1 \le q \le \hat q)$ and all residuals 
$\delta_k-\sum_{q=1}^{\hat q}\inprod{\C_k^q}{\widetilde{\V}^q} \ (k=1,2)$
are zero.
 In particular, $\widetilde{\V}^q=\O$ for every $q$, and hence $\mbox{rank}(\widetilde{\V}^q)=0 \le 1$ for every $q$. This proves the result.
\qed
}



\section{Examples}

\label{section:example}

This section presents two examples illustrating the main results. The first example is a separable homogeneous QCQP 
and shows when Assumption~(A) of Theorem~\ref{theorem:constLimited} holds. The second example illustrates 
how Theorem~\ref{theorem:main0} can be used to construct a separable QCQP with an exact SDPR 
by combining heterogeneous sub-QCQPs from different classes, including convex QCQPs, QCQPs 
characterized by sign-pattern conditions, and separable homogeneous QCQPs. It demonstrates a systematic integration within a unified framework and 
 shows how the exactness of the overall SDPR can be derived from the exactness of the individual subproblems.

\examp (A slight modification of \cite[Example 4.1]{KIM2023}) 
\label{example:constLimited}
This example illustrates the role of Assumption (A) in Theorem~\ref{theorem:constLimited}. 
In general, whether Assumption (A) holds depends on the problem data,
including the objective and constraint functions.
To make this dependence explicit, 
we introduce a parameter $\alpha \in \Real$ into the constraints
and consider the following parametric separable homogeneous QCQP.
Let 
\begin{eqnarray*}
& & g^1_0(\v) = (v_1)^2, \  g^1_1(\v) = (v_2)^2, \  
g^1_2(\v) =(v_1-\alpha v_2)(v_1-4v_2), \\ 
& & \hspace{26.5mm} g^1_3(\v) =  -(v_1-2v_2)(v_1-3v_2) \ \mbox{for every } \v^1 \in \Real^2,  \\
& & g^2_0(w) = -w^2, \ g^2_1(w) = g^2_2(w) \equiv 0, \  g^2_3(w) = w^2, \  \mbox{for every } w \in \Real, \  
 \bdelta = (1,0,0), 
\end{eqnarray*}
where $\alpha \in [0,4]$ is a parameter, which will be specified later. 
We consider the QCQP 
\begin{eqnarray}
\zeta(\bdelta) & = & \inf \left\{ g^1_0(\v) + g^2_0(w) : 
\begin{array}{l}
\v \in \Real^2, \ w \in \Real^1, \
g^1_1(\v) + g^2_1(w) = \delta_1, \\
g^1_2(\v) + g^2_2(w) \leq \delta_2, \  g^1_3(\v) + g^2_3(w) \leq \delta_3
\end{array}
\right\} \nonumber \\
& = & \inf \left\{ v_1^2 - w^2 : 
\begin{array}{l}
\v \in \Real^2, \ w \in \Real^1, \
v_2^2 = 1, \\
(v_1-\alpha v_2)(v_1-4v_2) \leq 0, \\ 
-(v_1-2v_2)(v_1-3v_2) + w^2 \leq 0
\end{array}
\right\} \label{eq:QCQPex1}\\
& = & \inf \left\{ v_1^2 - w^2 : 
\begin{array}{l}
\v \in \Real^2, \ w \in \Real^1, \ 
\alpha \leq v_1 \leq 4, \\ 
\displaystyle v_1 \leq \frac{5-\sqrt{1+4w^2}}{2} \ \mbox{or } \frac{5+\sqrt{1+4w^2}}{2} \leq v_1
\end{array}
\right\} \nonumber 
\end{eqnarray}
Here $v_2^2 = 1$ in the constraint of QCQP~\eqref{eq:QCQPex1} implies either $v_2 = 1$ or $v_2 = -1$. We have 
assumed $v_2 \geq 0$ without loss of generality since $g^1_k(\v) = g^1_k(-\v)$ $(0 \leq k \leq 3)$; hence $v_2=1$. 
This QCQP~\eqref{eq:QCQPex1} 
can be formulated as a special case of the separable homogeneous QCQP~\eqref{eq:QCQP3} with 
\begin{eqnarray*}
& & \hat{q} = 2,  \ m=3, \\
& & \C^1_0 = \begin{pmatrix} 1 & 0 \\ 0 & 0 \end{pmatrix}, \ 
\C^1_1 = \begin{pmatrix} 0 & 0 \\ 0 & 1 \end{pmatrix}, \
\C^1_2 = \begin{pmatrix} 1 & -\frac{4+\alpha}{2} \\ -\frac{4+\alpha}{2} & 4\alpha \end{pmatrix}, \
\C^1_3= \begin{pmatrix} -1 & \frac{5}{2} \\ \frac{5}{2} & -6 \end{pmatrix}, 
\\
& & \C^2_0 = -1, \C^2_1 = \C^2_2 = 0, \ \C^2_3 = 1, \ \mbox{`}\trianglelefteq_1\mbox{'} =  \mbox{`}=\mbox{'}, \   
\mbox{`}\trianglelefteq_2\mbox{'}=\mbox{`}\leq\mbox{'}, \mbox{`}\trianglelefteq_3\mbox{'} =  \mbox{`}\le\mbox{'}.
\end{eqnarray*}
We denote each feasible solution $(\V^1,\V^2) \in \SymMat^2_+ \times \SymMat^1_+$ 
of SDPR of QCQP~\eqref{eq:QCQP3} as $(\V,\W)$.  
The optimal solutions $(\tilde{\v},\tilde{\w})$ of QCQP~\eqref{eq:QCQPex1} and $(\widetilde{\V},\widetilde{\W})$ of 
its SDPR are summarized in Table 1. We note that $\widetilde{\W} = \tilde{w}^2$ holds. 
In particular, we compare the optimal  solution $(\widetilde{\V},\widetilde{\W})$ and optimal value of the SDPR 
for different values of $\alpha$.


\begin{table}[h!]
\begin{center}
\begin{tabular}{|c||c|c||c|c|c|c|c|c|}
\hline
              & \multicolumn{2}{|c||}{QCQP~\eqref{eq:QCQPex1}} &  \multicolumn{3}{|c|}{SDPR of QCQP~\eqref{eq:QCQPex1}} \\
 & Opt.sol. $(\tilde{\v},\tilde{w})$ & Opt.val. & Opt.val. & Opt.sol. $(\widetilde{\V},\widetilde{\W})$ & rank$(\widetilde{\V})$ \\ 
 \hline                
 $0 \le \alpha < 2$ & $\tilde{\v} = (\alpha,1), \tilde{w}\ne 0$ & $ \alpha^2-\tilde{w}^2$ &$\alpha^2-\tilde{w}^2$ & $\widetilde{\V} = \begin{pmatrix} \alpha^2 & \alpha \\ \alpha & 1 \end{pmatrix},  \widetilde{\W}  \ne \O$  & 1 \\
\hline                
 $\alpha =2 $ & $\tilde{\v} = (2,1), \tilde{w}=0$ & $4$ & $4$ & $\widetilde{\V} = \begin{pmatrix} 4 & 2 \\ 2 & 1 \end{pmatrix},  \widetilde{\W}= \O$ & 1 \\
\hline                
$2 < \alpha < 3 $ & $\tilde{\v} = (3,1), \tilde{w}=0$ &$9$ &$\frac{14\alpha-24}{\alpha-1}$ & $\widetilde{\V} = \begin{pmatrix} \frac{14\alpha-24}{\alpha-1} & \frac{4\alpha-6}{\alpha-1} \\ \frac{4\alpha-6}{\alpha-1} & 1 \end{pmatrix},\widetilde{\W}= \O$ & 2 \\
\hline 
$\alpha =3$ & $\tilde{\v} = (3,1), \tilde{w}=0$ & $9$ & $9$ & $\widetilde{\V} = \begin{pmatrix} 9 & 3 \\ 3 & 1 \end{pmatrix}, \widetilde{\W}= \O$ & 1 \\
\hline                           
$3 < \alpha \le 4 $ & $\tilde{\v} = (\alpha,1), \tilde{w}\ne 0$ &$\alpha^2-\tilde{w}^2 $ &$\alpha^2$ & $\widetilde{\V}=\begin{pmatrix} \alpha^2 & \alpha \\ \alpha & 1 \end{pmatrix},\widetilde{\W}\ne \O$ & 1\\
\hline       
\end{tabular}	
\end{center}
\caption{
Optimal solutions of QCQP (16) and its SDPR.
Here $(\tilde \v, \tilde w)$ denotes an optimal solution of QCQP~\eqref{eq:QCQPex1},
and $(\widetilde{\V}, \widetilde{\W})$ an optimal solution of its SDPR.
}
\end{table}


As shown in Table 1, the behavior of the SDPR  depends strongly on $\alpha$.
If $\alpha \in [0,2) \cup (3,4]$, then Assumption (A) of Theorem~\ref{theorem:constLimited} is satisfied, and 
hence the SDPR admits a rank-one optimal solution. 
In this case, the SDPR is exact. 
On the other hand, if $\alpha \in [2,3]$, Assumption (A) is violated. More precisely, the total 
number of nonzeros among 
\begin{eqnarray*}
\widetilde{\V}, \ 
\widetilde{\W}, \ \delta_k - \inprod{\C^1_k}{\widetilde{\V}} -\inprod{\C^2_k}{\widetilde{\W}} \ (k=1,2,3);  
\end{eqnarray*}
is less than $m-1=2$;  only $\widetilde{\V}$ is nonzero among these quantities. 
Particularly, if $\alpha \in (2,3)$,  
$\widetilde{\V}$ has rank two, and hence the SDPR is not exact. 
In the case where $\alpha = 2$ or $3$, we have rank$(\widetilde{\V}) = 1$
and rank$(\widetilde{\W}) = 0$. Although Assumption (A) is violated, the SDPR remains exact.
This shows that Assumption (A) is a sufficient but  not necessary  
for exactness.


The loss of exactness for $\alpha \in (2,3)$ can be understood
from the different behaviors of the feasible regions of QCQP~\eqref{eq:QCQPex1}
and its SDPR.
At $\alpha = 2$, the feasible region of QCQP~\eqref{eq:QCQPex1} undergoes a
discontinuous structural change. 
In contrast, the feasible region of the SDPR varies continuously 
with respect to $\alpha$. As a result, in the range $\alpha \in (2,3)$, 
the SDP relaxation no longer captures the geometry of the original 
problem exactly, leading to a rank-two optimal solution.
This example demonstrates that Assumption (A) is not merely technical,
but is essential for guaranteeing exactness. It also shows that exactness 
may fail even for a small perturbation of constraints.

\eexamp


\examp 
\label{example:main0} 
We construct a separable QCQP with exact SDPR by combining three heterogeneous sub-QCQPs:
 a convex QCQP for $p=1$, a QCQP characterized by sign-pattern conditions for $p=2$, and a separable homogeneous QCQP for $p=3$. 
More precisely, we connect two sub-QCQPs of the form~\eqref{eq:QCQPp} (for $p=1,2$) and 
one separable homogeneous sub-QCQP of the form~\eqref{eq:QCQP3} with $\hat q=3$ for $p=3$, and obtain the QCQP
\begin{eqnarray*}
\ \hspace{-15mm} \zeta(\bgamma)  = \inf \left\{ \sum_{p=1}^3 f^p_0(\u^p):
\begin{array}{l}
\u^p \in \Real^{n^p+1}, \ 
\displaystyle \sum_{p=1}^3 f^p_k(\u^p) \trianglelefteq_k \gamma_k\
(p=1,2,3,1 \leq k \leq 7)
\end{array}
\right\} 
\end{eqnarray*}
\begin{eqnarray}
= \inf  \left\{ \sum_{p=1}^2 f^p_0(\u^p)+ \sum_{q=1}^3 \inprod{\C^q_0}{\v^q(\v^q)^T} :
\begin{array}{l}
\u^p \in \Real^{n^p}, \v^q \in \Real^{\bar{n}^q} \
(p=1,2,q=1,2,3),\\
\displaystyle \sum_{p=1}^2 f^p_k(\u^p)+ \sum_{q=1}^3 \inprod{\C^q_k}{\v^q(\v^q)^T} \trianglelefteq_k \gamma_k\\
(1 \leq k \leq 7)
\end{array}
\right\}, \label{eq:QCQPexample}
\end{eqnarray}
and its SDPR
\begin{eqnarray}
\eta(\bgamma)  = \inf \left\{ \sum_{p=1}^2 \inprod{\B^p_0}{\X^p} + \sum_{q=1}^3 \inprod{\C^q_0}{\V^q} : 
\begin{array}{l}
\X^p \in \SymMat^{n^p+1}_+,\V^q \in \SymMat^{\bar{n}^q}_+,\\
X^p_{n^p+1n^p+1}=1\ (p=1,2, q=1,2,3),\\
\displaystyle \sum_{p=1}^2 \inprod{\B^p_k}{\X^p} + \sum_{q=1}^3 \inprod{\C^q_k}{\V^q} \trianglelefteq_k \gamma_k\\
(1 \leq k \leq 7)
\end{array}
\right\},  \label{eq:SDPRexample}
\end{eqnarray}
where 
\begin{eqnarray}
& & f^p_k(\u^p) = \inprod{\B^p_k}{\begin{pmatrix}\u^p\\1\end{pmatrix}\begin{pmatrix}\u^p\\1\end{pmatrix}^T} \ 
\ \mbox{for every } \u^p \in \Real^{n^p} \ (p=1,2, \ 0\leq k \leq 7), 
\nonumber \\
& & f^3_k(\u^3) = \sum_{q=1}^3\inprod{\C^q_k}{\v^q(\v^q)^T}  \ \mbox{for every } \u^3 = (\v^1,\v^2,\v^3) \in \Real^{n^3} \ (q=1,2,3,0\le k \le 7),
\nonumber \\
& & \B^p_k \in \SymMat^{n^p+1}  \ (p=1,2), \ \C^q_k \in \SymMat^{\bar{n}^q} \ (q=1,2,3), \ n^3 = \bar{n}^1+\bar{n}^2+\bar{n}^3 \ 
(0\le k \le 7).  \nonumber  
\end{eqnarray}
We assume that
\begin{eqnarray}
& & f^1_k \ \mbox{is convex} \ (0 \leq k \leq 7), f_1^1 \equiv f_2^1 \equiv 0 \ \mbox{(equivalently, }  \B_1^1=\B_2^1=\O\mbox{)}, \label{eq:convex2} \\
& & \mbox{all off-diagonal elements of } \B^2_k \ (0 \leq k \leq 7) \ \mbox{are nonpositive, } \B^2_1=\B^2_2 = \O, \label{eq:signPattern2} \\
& & \C^2_1 = \O, \C^3_1 = \O, 
\trianglelefteq_1 = \mbox{`}=\mbox{'},\gamma_1 \ne 0, \label{eq:const1} \\
& & \C^1_2,\C^3_2 \mbox{ are negative semidefinite},
\trianglelefteq_2 = \mbox{`}\ge\mbox{'},\gamma_2 > 0, \label{eq:const2}, \\
& & \B^1_3,\B^2_3,\C^1_3,\C^2_2 \mbox{ are positive semidefinite},\trianglelefteq_3 = \mbox{`}\le\mbox{'},\gamma_3 < 0, \label{eq:const3} \\
& & \alpha > 0, \ \C^q_5 = \alpha \C^q_4 \ (q=1,2,3), \label{eq:const4} \\
& & \C^q_k = \O \ (1 \leq q \leq 3,k=6,7), \label{eq:const5}\\ 
& & \trianglelefteq_k = \mbox{`}\leq\mbox{'} \ (3 \leq k \leq 7).  \label{eq:trianglelefteq1} 
\end{eqnarray}
The conditions imposed on the data matrices and constraints are summarized
in Table~\ref{Table2}, where  the structural properties of each sub-QCQP are included.
Table~2 shows how the three different classes of QCQPs are combined within a unified separable structure.


To apply Theorem~\ref{theorem:main0}, 
suppose that SDPR~\eqref{eq:SDPRexample} admits an optimal solution, as required by Assumption~(II), and denote it by
$(\widetilde{\X}^1,\widetilde{\X}^2,\widetilde{\V}^1,\widetilde{\V}^2,\widetilde{\V}^3)$.
 Let 
\begin{eqnarray*}
& & \delta_k^1 = \gamma_k - \inprod{\B^2_k}{\widetilde{\X}^2} - \sum_{q=1}^3 \inprod{\C^q_k}{\widetilde{\V}^q}, \
\delta_k^2 = \gamma_k - \inprod{\B1_k}{\widetilde{\X}^1} - \sum_{q=1}^3 \inprod{\C^q_k}{\widetilde{\V}^q}, \\
& & \delta_k^3 = \gamma_k - \sum_{p=1}^2 \inprod{\B^p_k}{\widetilde{\X}^p} \  \ (1 \leq k \leq 7).
\end{eqnarray*}
As shown in the proof of Theorem~\ref{theorem:main0}, $\widetilde{\X}^1$ is optimal for the SDPR of the first sub-QCQP with right-hand side $\bdelta^1$, $\widetilde{\X}^2$ is 
optimal for the SDPR of the second sub-QCQP with right-hand side $\bdelta^2$, and $(\widetilde{\V}^1,\widetilde{\V}^2,\widetilde{\V}^3)$ is optimal 
for SDPR~\eqref{eq:SDPR3} with right-hand side $\bdelta^3$.
%

For $p=1$, the first two constraints are variable-free because $\B_1^1=\B_2^1=\O$. 
After removing these constraints, the remaining sub-QCQP is convex by~\eqref{eq:convex2} and~\eqref{eq:trianglelefteq1}, 
and hence its SDPR is exact by Theorem~\ref{theorem:convex}. For $p=2$, 
the first two constraints are again variable-free because $\B_1^2=\B_2^2=\O$. 
After removing them, the remaining sub-QCQP satisfies the sign-pattern condition 
in Corollary~\ref{corollary:signDefinite} by~\eqref{eq:signPattern2} and~\eqref{eq:trianglelefteq1}, and therefore its SDPR is exact.



For $p=3$, sub-QCQP~\eqref{eq:QCQP3} is a separable homogeneous QCQP. We verify that its SDPR satisfies the assumptions of 
Theorem~\ref{theorem:constLimited}. 
By~\eqref{eq:const5}, the $6$th and $7$th constraints 
(the inequality constraints corresponding to $k = 6,7$) 
can be removed without affecting the feasible region of sub-SDPR~\eqref{eq:SDPR3}. 
By~\eqref{eq:const4}, the $4$th and $5$th inequality constraints,  
\[
\sum_{q=1}^3 \langle \C_4^q, \V^q \rangle \le \delta^3_4 \mbox{ and }  \sum_{q=1}^3 \langle \C_5^q, \V^q \rangle \le \delta^3_5
\]
have proportional coefficient matrices
with $\C_5^q=\alpha \C_4^q$ and $\alpha>0$. 
 Therefore one of them is redundant: the 4th constraint may be removed if $\delta_4^3 \le \delta_5^3/\alpha$, and the 5th otherwise. 
 Thus the sub-SDPR~\eqref{eq:SDPR3} reduces to a problem with at most four essential constraints.


Moreover, 
 every optimal solution satisfies $\V^1 \neq \O$, $\V^2 \neq \O$, and $\V^3 \neq \O$. Indeed, \eqref{eq:const1} forces $\V^1 \neq \O$, 
since $\gamma_1 \neq 0$ while $\C_1^2=\C_1^3=\O$; \eqref{eq:const2} forces $\V^2 \neq \O$, since $\gamma_2>0$ 
while $\C_2^1$ and $\C_2^3$ are negative semidefinite; and \eqref{eq:const3} forces $\V^3 \neq \O$, since $\gamma_3<0$ 
while $\C_3^1$ and $\C_3^2$ are positive semidefinite. 
Hence, after the above preprocessing, Assumption~(A) of Theorem~\ref{theorem:constLimited} is satisfied.



However, the right-hand-side values $\delta_4^3$ and $\delta_5^3$ depend on the optimal solution of 
the  full SDPR~\eqref{eq:SDPRexample}, in particular on $\widetilde{\X}^1$ and $\widetilde{\X}^2$. Therefore, 
it is not known a priori which of the 4th or 5th constraints can be removed in the analysis of the sub-SDPR~\eqref{eq:SDPR3}. 
Moreover, these two constraints correspond to different coupling constraints in SDPR~\eqref{eq:SDPRexample}, 
so both are essential in SDPR~\eqref{eq:SDPRexample}. 
 In particular, 
removing either of them would change the feasible region of SDPR~\eqref{eq:SDPRexample}. The elimination of one of these 
constraints is carried out solely within the analysis of the sub-SDPR~\eqref{eq:SDPR3},
 where it serves as a technical simplification, and should not be interpreted as indicating redundancy 
 in the original problem.


Therefore, each of the three sub-QCQPs satisfies Assumption~(I) of Theorem~\ref{theorem:main0}. 
We have shown that the SDPR of each sub-QCQP of QCQP~\eqref{eq:QCQPexample} is exact:
the case $p=1$ 
follows from convexity, the case $p=2$ from Corollary~\ref{corollary:signDefinite}, and the case $p=3$ 
from Theorem~\ref{theorem:constLimited}  after the above preprocessing. 
Therefore, Assumption (I) of Theorem~\ref{theorem:main0} is fulfilled for QCQP~\eqref{eq:QCQPexample}
and its SDPR~\eqref{eq:SDPRexample}. 
By applying Theorem~\ref{theorem:main0}, 
we conclude that SDPR~\eqref{eq:SDPRexample} is an exact relaxation of QCQP~\eqref{eq:QCQPexample}, 
{\it i.e.},  both problems attain the same optimal value, 
whenever SDPR~\eqref{eq:SDPRexample} has an optimal solution (Assumption (II)). 
\eexamp

\begin{table}[h!] 
\begin{center}
\begin{tabular}{|c|c|c||c|c|c||c|c|c|} 
\hline
           &  $p=1$: Convex                 & $p=2$:  Sign Pattern                 & \multicolumn{3}{|c||}{$p=3$: Homogeneous Separable} & & \\
           & Sub-QCQP~\eqref{eq:QCQPp} & Sub-QCQP~\eqref{eq:QCQPp}  &  \multicolumn{3}{|c||}{Sub-QCQP~\eqref{eq:QCQP3}} & & \\ 
 $k$     &  See Section~\ref{section:convexity}    & See Section~\ref{section:signPattern} &  \multicolumn{3}{|c||}{See Section~\ref{section:separable}} & $\trianglelefteq_k$ & $\gamma_k$ \\ 
\hline 
  $0$         &  $f^1_0$:convex   &  $\B^2_0$:Off-diag$\ominus$   &  $\forall\C^1_0$              &  $\forall\C^2_0$               & $\forall\C^3_0$    &  & \\ \hline 
  $1 $         &  $\B^1_1=\O$  & $\B^2_1=\O$   &  $\forall\C^1_1$ & $\C^2_1=\O$                      & $\C^3_1=\O$          & $=$ & $\ne 0$\\
\hline 
  $2$          &  $\B^1_2=\O$  & $\B^2_2=\O$   &  $\C^1_2$:$\ominus$                    & $\forall\C^2_2$    & $\C^3_2$:$\ominus$  
  &  $\geq$  & $+ $ \\
 \hline
   $3$ &  $f^1_3$:convex                           &  $\B^2_3$:Off-diag$\ominus$               &                               &                    & &  &\\
          & $\B^1_3$:$\oplus$  & $\B^2_3$:$\oplus$   &  $\C^1_3$:$\oplus$     &  $\C^2_3$:$\oplus$   & $\forall\C^3_3$  & $\leq$ & $-$ \\ 
\hline   
$4$  &  $f^1_4$:convex             & $\B^2_4$:Off-diag$\ominus$ & $\forall\C^1_4$              &  $\forall\C^2_4$               & $\forall\C^3_4$             &  $\le$  & $\forall$ \\
\hline
$5$      & $f^1_5$:convex              & $\B^2_5$:Off-diag$\ominus$ &  $\C^1_{5} = \alpha\C^1_4$ &  $\C^2_{5} = \alpha\C^2_4$ & $\C^3_{5} = \alpha\C^3_4$ & $\le$  & $\forall$ \\
\hline
$6$ &   $f^1_6$:convex             & $\B^2_6$:Off-diag$\ominus$ & $\C^1_6=\O$  &  $\C^2_6=\O$ & $\C^3_6=\O$ &  $\le$  & $\forall$ \\
\hline
$7$ &   $f^1_7$:convex              & $\B^2_7$:Off-diag$\ominus$ & $\C^1_7=\O$  &  $\C^2_7=\O$  & $\C^3_7=\O$&  $\le$  & $\forall$ \\
\hline                
\end{tabular}	
\end{center}
\caption{
Structural conditions for the three sub-QCQPs.
The table summarizes the assumptions (19)–(26), indicating
how each constraint contributes to the convexity,
sign-pattern condition, or separable homogeneous structure.
$\oplus$ : Positive semidefinite.  $\ominus$ : Negative semidefinite.  
Off-diag$\ominus$ : Off-diagonal nonpositive. 
$\alpha \geq 0$. 
}    \label{Table2}
\end{table}


\section{Concluding remarks}

\label{section:ConcludingRemarks}

In this paper, we have shown that exactness of the SDPR for a class of 
separable QCQPs  is preserved under a 
horizontal connection of sub-QCQPs whose 
SDPRs are exact. Theorem~\ref{theorem:main0} provides a simple framework for constructing separable 
QCQPs by coupling sub-QCQPs through their right-hand-side constraint parameters. 
Rather than introducing a fundamentally new structural class, the framework clarifies 
how exactness results for existing QCQP classes can be systematically combined through separable structure.
In addition, Theorem~\ref{theorem:constLimited} presents a rank-based exactness condition 
for separable homogeneous QCQPs, 
which may  be viewed as a refinement of earlier results based on the bound 
$m \le \hat{q} + 1$.
These results illustrate how exactness can be ensured in structured separable QCQPs.


Section~\ref{section:subQCQPs} presents three classes of QCQPs satisfying Assumption (I) of
Theorem~\ref{theorem:main0}, based on convexity, sign-pattern and graph-structural conditions,
and separability with a limited number of constraints. Although the
mechanisms guaranteeing exactness differ across these classes, they can all
serve as sub-QCQPs in the horizontal connection of Theorem~\ref{theorem:main0}, leading to new
separable QCQPs with exact SDPRs. 


This horizontal 
connection is complementary to the vertical extension framework of Kojima, Kim, and Arima \cite{KOJIMA2025},
which is particularly relevant to QCQPs arising from classes (d) and (e). 
Once the exactness of an SDPR is established for a given QCQP, the
vertical extension framework  can be applied to
further enlarge the class of QCQPs by adding quadratic inequality
constraints, while preserving exactness under suitable NIQC assumptions. 
Together, these horizontal and vertical extensions suggest a systematic way to construct structured QCQPs 
whose SDPR exactness is guaranteed by construction.

%
%


\end{document}